# Spinor Bishop Equations of Curves in Euclidean 3-Space

Doğan ÜNAL[1]    İlim KİŞİ[2]    Murat TOSUN[1]

1. Department of Mathematics, Sakarya University, Sakarya 54187, Turkey

2. Department of Mathematics, Kocaeli University, Kocaeli 41100, Turkey

## Abstract

In this paper, we study spinor Bishop equations of curves in $\mathbb{E}^3$. We research the spinor formulations of curves according to Bishop frames in $\mathbb{E}^3$. Also, the relation between spinor formulations of Bishop frames and Frenet frame are expressed.

**AMS 2010:** 55A04, 15A66.

## 1 Introduction

Spinors in general were discovered by Elie Cartan in 1913, [2]. Later, spinors were adopted by quantum mechanics in order to study the properties of the intrinsic angular momentum of the electron and other fermions. Today, spinors enjoy a wide range of physics applications. In mathematics, particularly in differential geometry and global analysis, spinors have since found broad applications to algebraic and differential topology, symplectic geometry, gauge theory and complex algebric geometry, [4,7].

The Frenet frame is constructed for the curve of 3-time continuously differentiable non-degenerate curves. But, curvature may vanish at some points on the curve. That is, second derivative of the curve may be zero. In this situation, we need an alternative frame in $\mathbb{E}^3$. Therefore, in [1], Bishop defined a new frame for a curve and he called it Bishop frame (parallel transport frame) which is well defined even if the curve has vanishing second derivative in Euclidean 3-space. The advantages of Bishop frame and comparable Bishop frame with the Frenet frame in Euclidean 3-space was given by Bishop [1] and Hanson [6].

In [2,13], the triads of mutually orthogonal unit vectors were expressed in terms of a single vector with two complex components, called a spinor. In the light of the existing studies, the Frenet equations reduce to a single spinor equation, equivalent to the three usual vector equations, is a consequence of the relationship between spinors and orthogonal triads of vectors. The aim of this paper to show that the Bishop equations can be expressed a single eqation for a vector with two complex components.

## 2 Preliminaries

The Euclidean 3-space provided with the standard flat metric given by

$$\langle\ ,\ \rangle = dx_1^2 + dx_2^2 + dx_3^2 \tag{1}$$

where $(x_1, x_2, x_3)$ is a rectangular coordinate system of $E^3$. Recall that the norm of an arbitrary vector $X \in \mathbb{E}^3$ is given by $\|X\| = \sqrt{\langle X, X \rangle}$. Let $\alpha: I \subset \mathbb{R} \longrightarrow \mathbb{E}^3$ be arbitrary curve in the Euclidean space $\mathbb{E}^3$. The curve $\alpha$ is said to be of unit speed (parametrized by arc length function $s$) if $\langle \frac{d\alpha}{ds}, \frac{d\alpha}{ds} \rangle = 1$, the derivatives of the Frenet frame (Frenet-Serret formula),

$$\frac{dT}{ds} = \kappa N$$

$$\frac{dN}{ds} = -\kappa T + \tau B \tag{2}$$

$$\frac{dB}{ds} = -\tau N$$

where $\{T, N, B\}$ is Frenet frame of $\alpha$ and $\kappa, \tau$ are the curvature and the torsion of the curve $\alpha$, respectively, [3].

The Bishop frame or parallel transport frame is an alternative approach to defining a moving frame that is well-define even when the curve has vanishing second derivative. One can express parallel transport of an orthonormal frame along a curve simply by parallel transporting each component of the frame. The tangent vector and any convenient arbitrary basis for the remainder of the frame are used. Therefore, the type-1 Bishop (frame) formulas is expressed as

$$\frac{dT}{ds} = k_1 N_1 + k_2 N_2$$

$$\frac{dN_1}{ds} = -k_1 T \tag{3}$$

$$\frac{dN_2}{ds} = -k_2 T$$

Here, we shall call the set $\{T, N_1, N_2\}$ as type-1 Bishop frame and $k_1$ and $k_2$ are called the first and second Bishop curvature, respectively, [1].

The relation between Frenet frame and type-1 Bishop frame is given as follows.

$$T = T$$
$$N = cos\theta N_1 + sin\theta N_2 \qquad (4)$$
$$B = -sin\theta N_1 + cos\theta N_2$$

where $\theta(s) = arctan\left(\frac{k_2}{k_1}\right)$, $\tau(s) = \left(\frac{d\theta(s)}{ds}\right)$ and $\kappa(s) = \sqrt{k_1^2 + k_2^2}$. Here type-1 Bishop curvatures are defined by

$$k_1 = \kappa cos\theta \qquad (5)$$
$$k_2 = \kappa sin\theta.$$

Furthermore, the type-2 Bishop (frame) formulas of the curve $\alpha$ is defined by

$$\frac{d\zeta_1}{ds} = -\varepsilon_1 B$$
$$\frac{d\zeta_2}{ds} = -\varepsilon_2 B \qquad (6)$$
$$\frac{dB}{ds} = \varepsilon_1 \zeta_1 + \varepsilon_2 \zeta_2$$

Here, we call the set $\{\zeta_1, \zeta_2, B\}$ as type-2 Bishop frame, $\varepsilon_1$ and $\varepsilon_2$ are called type-2 Bishop curvatures, [9].

The relation between Frenet and type-2 Bishop frames can be expressed.

$$T = sin\theta \zeta_1 - cos\theta \zeta_2$$
$$N = cos\theta \zeta_1 + sin\theta \zeta_2 \qquad (7)$$
$$B = B$$

Here the type-2 Bishop curvatures are defined by

$$\varepsilon_1 = -\tau cos\theta \qquad (8)$$
$$\varepsilon_2 = -\tau sin\theta$$

It can be also deduced as

$$\theta'(s) = \kappa = \frac{(\varepsilon_2/\varepsilon_1)'}{1 + (\varepsilon_2/\varepsilon_1)^2} \qquad (9)$$

The frame $\{\zeta_1, \zeta_2, B\}$ is properly oriented and $\tau$ and $\theta(s) = \int_0^s \kappa(s)ds$ are polar coordinates for the curve $\alpha(s)$.

# 3 Spinors and Orthonormal Bases

The group of rotation about the origin denoted by $SO(3)$ in $\mathbb{R}^3$ is homomorphic to the group of unitary complex 2x2 matrices with unit determinant denoted by $SU(2)$.

Therefore, there is a two-to-one homomorphism of $SU(2)$ on to $SO(3)$. Whereas the elements of $SO(3)$ act the vectors with three real components (the points of $\mathbb{R}^3$), the elements of $SU(2)$ act on vectors with two complex components which are called spinors, [5,11].

We can define a spinor

$$\Psi = \begin{pmatrix} \Psi_1 \\ \Psi_2 \end{pmatrix} \tag{10}$$

by means of three vectors $a, b, c \in \mathbb{R}^3$ such that

$$a + ib = \Psi^t \sigma \Psi, \quad c = -\widehat{\Psi}^t \sigma \Psi \tag{11}$$

where $\sigma = (\sigma_1, \sigma_2, \sigma_3)$ is a vector whose cartesian components are the complex symmetric 2x2 matrices

$$\sigma_1 = \begin{pmatrix} 1 & 0 \\ 0 & -1 \end{pmatrix}, \quad \sigma_2 = \begin{pmatrix} i & 0 \\ 0 & i \end{pmatrix}, \quad \sigma_3 = \begin{pmatrix} 0 & -1 \\ -1 & 0 \end{pmatrix} \tag{12}$$

Let $\widehat{\Psi}$ be the mate (or conjugate) of $\Psi$ and $\overline{\Psi}$ be complex conjugation of $\Psi$, [10]. Therefore,

$$\widehat{\Psi} = -\begin{pmatrix} 0 & 1 \\ -1 & 0 \end{pmatrix} \overline{\Psi} = \begin{pmatrix} 0 & 1 \\ -1 & 0 \end{pmatrix} \begin{pmatrix} \overline{\Psi}_1 \\ \overline{\Psi}_2 \end{pmatrix} = \begin{pmatrix} -\overline{\Psi}_2 \\ \overline{\Psi}_1 \end{pmatrix} \tag{13}$$

Taking $a + ib = (x_1, x_2, x_3)$, from the equation (11) and (12) we have

$$x_1 = \Psi^t \sigma_1 \Psi = \Psi_1^2 - \Psi_2^2, \quad x_2 = \Psi^t \sigma_2 \Psi = i(\Psi_1^2 + \Psi_2^2), \quad x_3 = \Psi^t \sigma_3 \Psi = -2\Psi_1 \Psi_2$$

where superscript $t$ denotes transposition. That is,

$$a + ib = \Psi^t \sigma \Psi = (\Psi_1^2 - \Psi_2^2, i(\Psi_1^2 + \Psi_2^2), -2\Psi_1 \Psi_2)$$

In the same manner we can see that

$$c = (c_1, c_2, c_3) = (\Psi_1 \overline{\Psi}_2 + \overline{\Psi}_1 \Psi_2, i(\Psi_1 \overline{\Psi}_2 - \overline{\Psi}_1 \Psi_2), |\Psi_1|^2 - |\Psi_2|^2)$$

Since the vector $a + ib \in \mathbb{C}^3$ is a isotropic vector, by means of an easy computation one find that $a, b$ and $c$ mutually orthogonal. Furthermore $|a| = |b| = |c| = \overline{\Psi}^t \Psi$ and $\langle a \wedge b, c \rangle = \det(a, b, c) > 0$.

On the contrary, if the vectors $a, b$ and $c$ mutually orthogonal vectors of same magnitude such that $\det(a, b, c) > 0$, then there exists a spinor, defined up to sign such that the equation (11) holds.

Under the condition stated above, for two arbitrary spinors $\phi$ and $\Psi$, there exist following equalities

$$\overline{\phi^t \sigma \Psi} = -\hat{\phi}^t \sigma \hat{\Psi}$$

$$\widehat{a\phi + b\Psi} = \bar{a}\hat{\phi} + \bar{b}\hat{\Psi} \qquad (14)$$

$$\hat{\hat{\Psi}} = -\Psi$$

where $a$ and $b$ are complex numbers [2,13].

The correspondence between the spinors and the orthogonal bases (given by the equation (11)) is two-to-one because the spinors $\Psi$ and $-\Psi$ correspond to the same ordered orthogonal bases $\{a, b, c\}$ with $|a| = |b| = |c|$ and $\langle a \wedge b, c \rangle > 0$.

In addition to that, the ordered triads $\{a, b, c\}, \{b, c, a\}, \{c, a, b\}$ correspond to different spinors. Since the matrices $\sigma$ (given by equation (12)) are symmetric, any pair of spinors $\phi$ and $\Psi$ satisfying

$$\phi^t \sigma \Psi = \Psi^t \sigma \phi \qquad (15)$$

The set $\{\hat{\Psi}, \Psi\}$ is linearly independent for the spinor $\Psi \neq 0$, [2,13].

## 4 Spinor Bishop Equations

In this section, the following similar way as spinor formulation of Frenet Frame in [6], the spinor formulation of Bishop frames of the space curve are investigated.

Let $\alpha: I \subset \mathbb{R} \to \mathbb{R}^3$ be any curve with the arc length parameterization and $T(s), N(s)$ and $B(s)$ be the unit tangent vector, unit normal vector and unit binormal vector of α, respectively. In this case, we know that the Frenet formulas is given by

$$\frac{dT}{ds} = \kappa N$$

$$\frac{dN}{ds} = -\kappa T + \tau B \qquad (16)$$

$$\frac{dB}{ds} = -\tau N$$

where $\kappa(s)$ and $\tau(s)$ are called the curvature and torsion of the curve $\alpha(s)$, respectively, [3].

According to the results concerned with the spinor (given by section 3) there are a spinor $\Psi$ such that

$$N + iB = \Psi^t \sigma \Psi \quad , T = -\hat{\Psi}^t \sigma \Psi \qquad (17)$$

with $\hat{\Psi}^t \Psi = 1$. Therefore, the spinor $\Psi$ represent the triad $\{N, B, T\}$ and the variations of this triad along the curve $\alpha$ must correspond to some expression for $\frac{d\Psi}{ds}$. That is, the Frenet equations are equivalent to the single spinor equation

$$\frac{d\Psi}{ds} = \frac{1}{2}\left(-i\tau\Psi + \kappa\hat{\Psi}\right) \qquad (18)$$

where $\kappa$ and $\tau$ denote the curvature and the torsion of the curve $\alpha$, respectively. The equation (18) is called spinor Frenet equation, [13].

The parallel transport frame or Bishop frame is based on the observation that, while $T(s)$ for a given curve model is unique, we may choose any convenient arbitrary basis. $(N_1(s), N_2(s))$ for the remainder of the frame, so long as it is in the normal plane perpendicular to $T(s)$ at each point. If the derivatives of $(N_1(s), N_2(s))$ depend only on $T(s)$ and not each other. We can make $N_1(s)$ and $N_2(s)$ vary smoothly throughout the path regardless of curvature.

Now, we consider type-1 and type-2 Bishop frames which is given by equations (3) and (6), respectively.

Let $\{T, N_1, N_2\}$ be type-1 Bishop frame of the unit space curve $\alpha$. From equation (3), we know that type-1 Bishop frame is as follows

$$\begin{aligned}\frac{dT}{ds} &= k_1 N_1 + k_2 N_2 \\ \frac{dN_1}{ds} &= -k_1 T \\ \frac{dN_2}{ds} &= -k_2 T\end{aligned} \qquad (19)$$

In this case, there are a spinor $\phi$ such that

$$N_1 + iN_2 = \phi^t \sigma \phi, \quad T = -\hat{\phi}^t \sigma \phi \qquad (20)$$

with $\bar{\phi}^t \phi = 1$ and the spinor $\phi$ denotes the triad $\{N_1, N_2, T\}$. In addition to that, the variations of the triad $\{N_1, N_2, T\}$ along the curve must correspond to the expression for $\frac{d\phi}{ds}$. Since $\{\phi, \hat{\phi}\}$ is a basis for the two component spinors ($\phi \neq 0$), there are two functions $f$ and $g$ such that

$$\frac{d\phi}{ds} = f\phi + g\hat{\phi} \qquad (21)$$

where the functions $f$ and $g$ are possibly complex valued functions.

Differentiating the first equation (20) with respect to $s$, we have

$$\frac{dN_1}{ds} + i\frac{dN_2}{ds} = \frac{d}{ds}(\phi^t \sigma \phi) = (\frac{d\phi}{ds})^t \sigma \phi + \phi^t \sigma \left(\frac{d\phi}{ds}\right) \qquad (22)$$

Substituting the equations (19), (20) and (21) in to (22) and after simplifying, one finds

$$-k_1 T - ik_2 T = 2 f (N_1 + iN_2) - 2 g (T) \qquad (23)$$

From the last equation it is obvious that

$$f = 0, \qquad g = \left(\frac{k_1 + i\, k_2}{2}\right) \qquad (24)$$

In this case, from equations (21) and (24) we give following theorem.

**Theorem 1:** Let the two-component spinors $\phi$ represents the triad $\{N_1, N_2, T\}$ of a curve parametrized by arc length. Then, the type-1 Bishop frame is equivalent to the single spinor equation.

$$\frac{d\phi}{ds} = \frac{1}{2}\, (k_1 + i\, k_2)\, \hat{\phi} \qquad (25)$$

where $k_1$ and $k_2$ are Bishop curvatures of the curve.

Now, we investigate the relation between the spinors $\Psi$ and $\phi$ represent the Frenet frame $\{N, B, T\}$ and the type-1 Bishop frame $\{N_1, N_2, T\}$, respectively.

Considering equation (4) we reach

$$N_1 + iN_2 = (N + \dot{\imath} B)(\cos\theta + i\sin\theta) \qquad (26)$$

From the equations (17) and (20) and the last equation, it easily seen that

$$\phi^2 \sigma \phi = e^{i\theta}(\Psi^t \sigma \Psi)$$
$$T = T \qquad (27)$$

Thus, we have proved the following theorem.

**Theorem 2:** The relations between the spinor formulations of type-1 Bishop frame and Frenet Frame are as follows.

$$\phi^t \sigma \phi = e^{i\theta}(\Psi^t \sigma \Psi)$$
$$T = T$$

where the spinors $\Psi$ and $\phi$ represent the Frenet frame $\{N, B, T\}$ and type-1 Bishop frame $\{N_1, N_2, T\}$, respectively.

We now apply to same argument again, with type-1 Bishop frame replaced by type-2 Bishop frame, to obtain the spinor formulation of type-2 Bishop frame.

From equation (11), there exists a spinor $\lambda$ such that

$$\zeta_1 + i\zeta_2 = \lambda^t \sigma \lambda, \quad B = -\hat{\lambda}^t \sigma \lambda \tag{28}$$

where $\overline{\lambda}^t \lambda = 1$ and the spinor $\lambda$ denotes the triad $\{\zeta_1, \zeta_2, B\}$.

Differentiating the equation of (28) and using type-2 Bishop equations (6) one finds that

$$-\varepsilon_1 B - i\varepsilon_2 B = 2f(\zeta_1 + i\zeta_2) - 2g(B)$$

which amounts to $f = 0$ and $g = \frac{\varepsilon_1 + \varepsilon_2}{2}$.

In this case, we obtain

$$\frac{d\lambda}{ds} = g\hat{\lambda} = \left(\frac{\varepsilon_1 + \varepsilon_2}{2}\right)\hat{\lambda} \tag{29}$$

such that $\{\lambda, \hat{\lambda}\}$ is a basis for the two component spinors. In addition to that; using the equation (7) easily seen that

$$\zeta_1 + i\zeta_2 = -ie^{i\theta}(T + iN)$$

and

$$\lambda^t \sigma \lambda = -ie^{i\theta}(\Psi^t \sigma \Psi)$$

$$B = B$$

where the spinors $\lambda$ and $\Psi$ represent the triads $\{\zeta_1, \zeta_2, B\}$ and $\{T, N, B\}$, respectively.

Thus, we have proved following theorems.

**Theorem 3:** Let the spinors $\lambda$ represent the triad $\{\zeta_1, \zeta_2, B\}$ of a unit speed regular curve. Then, the type-2 Bishop frame is equivalent to the single spinor equation.

$$\frac{d\lambda}{ds} = \frac{1}{2}(\varepsilon_1 + i\varepsilon_2)\hat{\lambda} \tag{30}$$

where $\varepsilon_1$ and $\varepsilon_2$ are Bishop curvatures of the curve.

**Theorem 4:** The relations between the spinor formulation of type-2 Bishop frame and Frenet frame are

$$\lambda^t \sigma \lambda = -ie^{i\theta}(\Psi^t \sigma \Psi)$$

$$B = B$$

where $\lambda$ and $\Psi$ are spinors and they represent the triad $\{\zeta_1, \zeta_2, B\}$ and $\{T, N, B\}$, respectively.